\documentstyle[12pt,twoside]{article}
\newtheorem{theorem}{\bf Theorem}[section]
\newtheorem{proposition}[theorem]{\bf Proposition}
\newtheorem{lemma}[theorem]{\bf Lemma}
\newtheorem{corollary}[theorem]{\bf Corollary}
\newtheorem{example}[theorem]{\bf Example}

\newtheorem{remark}[theorem]{\bf Remark}

\input{amssym}

\date{}
\begin{document}

\title{{\Large\bf Regularity and amenability of weighted Banach algebras and their second dual on locally compact groups}}

\author{{\normalsize\sc M. J. Mehdipour and A. Rejali\footnote{Corresponding author}}}
\maketitle

{\footnotesize  {\bf Abstract.} Let $\omega$ be a weight function on a locally compact group $G$ and let $M_*(G, \omega)$ be the subspace of $M(G, \omega)^*$ consisting of all functionals that
vanish at infinity. In this paper, we first investigate the Arens regularity of $M_*(G, \omega)^*$ and show that $M_*(G, \omega)^*$ is Arens regular if and only if $G$ is finite or $\Omega$ is zero cluster. This result is an answer to the question posed and it improves some well-known results. We also give necessary and sufficient criteria for the weight function spaces $Wap(G, 1/\omega)$ and $Ap(G, 1/\omega)$ to be equal to $C_b(G, 1/\omega)$. We prove that for non-compact group $G$, the Banach algebra $M_*(G, \omega)^*$ is Arens regular if and only if $Wap(G, 1/\omega)=C_b(G, 1/\omega)$. We then investigate amenability of $M_*(G, \omega)^*$ and prove that $M_*(G, \omega)^*$ is amenable and Arens regular if and only if $G$ is finite.}
{}
\section{\normalsize\bf Introduction}

Throughout this paper, $G$ denotes a Hausdorff locally compact group with the group algebra $L^1(G)$ and the measure algebra $M(G)$. A \emph{weight} on $G$ is a continuous function $\omega: G\rightarrow [1, \infty)$ such that
$\omega(e)=1$ and $$\omega(xy)\leq\omega(x)\;\omega(y)$$ for all $x, y\in G$, where $e$ is the identity element of $G$.
Let the function $\Omega: G\times G\rightarrow (0, 1]$ be defined as follows:
$$
\Omega(x, y)=\omega(xy)/\omega(x)\omega(y).
$$
Let us recall that a complex-valued function $F$ on $G\times G$ is called \emph{cluster} (respectively, \emph{zero cluster, positive cluster}) if for every pair of sequences $(x_n)_n$ and $(y_m)_m$ of distinct elements in $G$, we have
\begin{eqnarray}\label{positive cluster}
\lim_n\lim_mF(x_n, y_m)=\lim_m\lim_nF(x_n, y_m),
\end{eqnarray}
(respectively, both limits equal zero, positive) whenever both iterated limits exist.

Let $L^1(G, \omega)$ be the
space of all measurable functions $\phi$ on $G$ such that
$\omega\phi\in L^1(G)$.
Let also $M(G, \omega)$ be the Banach space of all complex regular Borel measures
$\mu$ on $G$ for which $$\omega\mu\in M(G).$$ It is well-known that $L^1(G,\omega)$ and $M(G, \omega)$ are Banach algebra and
$M(G, \omega)$ is the dual
space of $C_0(G, 1/\omega)$, the Banach space of all
complex-valued continuous functions $f$ on $G$ such that $f/\omega$ vanishes at infinity, see for example \cite{dl, r0}.

We say that $\lambda \in M(G, \omega)^*$ \emph{vanishes at infinity}
if for every $\varepsilon>0$, there exists a compact
subset $K$ of $G$, for which $|\langle \lambda, \mu \rangle|<\varepsilon$, where
$\mu \in M(G, \omega)$ with $|\mu|(K)=0$ and $\|\mu\|_\omega=1$.
We denote by $M_*(G, \omega)$ the subspace of $M(G, \omega)^*$ consisting of all functionals that
vanish at infinity. In the case where, $\omega(x)=1$ for all $x\in G$, we write the spaces
$$
M_*(G, \omega):=M_*(G).
$$
The space $M_*(G, \omega)$ is a norm closed subspace of $M(G, \omega)^*$ and so it is a $C^*-$algebra. Every element $f\in C_0(G, 1/\omega)$ may be regarded as an element in $M_*(G, \omega)$ by the pairing
$$
\langle f, \mu\rangle=\int_G f d\mu\quad\quad (M(G, \omega)).
$$
Then $C_0(G, 1/\omega)$ is a closed subspace of $M_*(G, \omega)$. Also, the space $M_*(G, \omega)$ is left introverted in $M(G, \omega)^*$. This let us to endow $M_*(G, \omega)^*$ with the first Arens product. Then $M_*(G, \omega)^*$ with this product becomes to a Banach algebra \cite{mm}.
For each $\phi\in L^1(G, \omega)$, let $\phi$ denote the functional in $M_*(G, \omega)^*$ defined by
$$
\langle\phi, \lambda\rangle:=\langle\lambda, \phi\rangle.
$$
for all $\lambda\in M(G, \omega)^*$. This duality defines a linear isometric embedding from $L^1(G, \omega)$ into $M_*(G, \omega)^*$. One can prove that $L^1(G, \omega)$ is a closed ideal in $M_*(G, \omega)^*$ and  $M_*(G, \omega)^*= L^1(G, \omega)$ if and only if $G$ is discrete \cite{mm}; see \cite{m} for the case $\omega=1$. Since $M(G, \omega)$ is a closed subspace of $M_*(G, \omega)^*$, an easy application of the Goldstine's theorem shows that if $\Phi\in M_*(G, \omega)^*$, then there exists a net $(\mu_\alpha)_\alpha$ in $M(G, \omega)$ such that $\mu_\alpha\rightarrow\Phi$ in the weak$^*$-topology of $M_*(G, \omega)^*$.

Let us recall that the first Arens product ``$\diamond$" on the second dual of a Banach algebra $\frak{A}$ is defined by
$$
\langle\Phi\diamond\Psi, f\rangle=\langle\Phi, \Psi f\rangle,
$$
in which
$$
\langle\Psi f, a\rangle =\langle\Psi, fa\rangle\quad\hbox{and}\quad\langle fa, b\rangle=\langle f, ab\rangle
$$
for all $\Phi, \Psi\in \frak{A}^{**}$, $f\in\frak{A}^*$ and $a, b \in\frak{A}$.
The Banach algebra $\frak{A}$ is called \emph{Arens regular} if for every $\Phi\in \frak{A}^{**}$ the mapping $\Psi\mapsto\Phi\diamond\Psi$ is weak$^*-$weak$^*$ continuous on $\frak{A}^{**}$.

Several authors have studied the Arens regularity of weighted group algebras. For example, Crow and Young \cite{cy} showed that there exists a weighted function $\omega$ on $G$ such that $L^1(G, \omega)$ is Arens regular if and only if $G$ is discrete and countable. The second author and Vishki \cite{rv} proved that $L^1(G, \omega)$ is Arens regular if and only if $G$ is finite or $G$ is discrete and $\Omega$ is zero cluster. They showed that $L^1(G, \omega)$ is amenable and Arens regular if and only if $G$ is finite; see also \cite{f}. These studies have continued for the other Banach algebras. See for example, \cite{br, r1} for the Arens regularity of weighted semigroup algebras and \cite{f1, g1, g2, g3, lu} for the Arens regularity of Fourier algebras. See also \cite{d, u, u1, u2}.

In this paper, we investigate the Arens regularity of $M_*(G, \omega)^*$ and the relation between it, the weighted function spaces  and amenability. In Section 2, we give an answer to the question presented in \cite{mm} and prove that $M_*(G, \omega)^*$ is Arens regular if and only if $G$ is finite or $\Omega$ is zero cluster. This result is an improvement of Theorem 2 of \cite{rv}. We also show that $M_*(G)^*$ is Arens regular if and only if there exists a weight function $\omega$ on $G$ such that $M_*(G, \omega)^*$ is $C^*-$algebra; or equivalently, $G$ is finite. In Section 3, we prove that $G$ is weight regular if and only if $G$ is a countable discrete group. For a normal subgroup $N$ of $G$, we show that if $G$ is weight regular, then $G/N$ is weight regular and $N$ is countable and open. Section 4 is devote to weighted function spaces $\hbox{Wap}(G, 1/\omega)$ and $\hbox{Ap}(G, 1/\omega)$. We give necessary and sufficient condition for these weighted  function spaces to be equal to $C_b(G, 1/\omega)$. For instance, we show that $\hbox{Wap}(G, 1/\omega)=C_b(G, 1/\omega)$ if and only if $G$ is compact or $\Omega$ is zero cluster. As a consequence of this result, we prove that  $M_*(G, \omega)^*$ is Arens regular if and only if $\hbox{Wap}(G, 1/\omega)=C_b(G, 1/\omega)$ , when $G$ is non-compact. In Section 5, we investigate amenability of $M_*(G, \omega)^*$ and prove that $M_*(G, \omega)^*$ is amenable if and only if $G$ is a discrete amenable group and $\omega^*$ is bounded. We also show that $M_*(G, \omega)^*$ is Arens regular and amenable if and only if $G$ is finite.

\section{\normalsize\bf Arens regularity of $M_*(G, \omega)^*$}

The following lemma is needed to prove our results.

\begin{lemma}\label{lem13} Let $\omega$ be a weight function on a locally compact group $G$. If $\Omega$ is zero cluster, then $G$ is discrete.
\end{lemma}
{Proof.} Suppose that $G$ is a non-discrete group. Let $\frak{U}$ be the family of all neighborhood of $e$ directed by upward inclusion, i.e.,
$$
U_1\geq U_2\Leftrightarrow U_1\subseteq U_2\quad\quad (U_1, U_2\in\frak{U}).
$$
Assume that $U\in\frak{U}$. Since $G$ is non-discrete, $U$ is infinite. So we can choose $x_U\in U$ such that $x_U\neq e$. Then the net $(x_U)_{U\in\frak{U}}\;$ of distinct points of $G$ converges to the identity element $e$. Indeed, if $W$ is a neighborhood of $e$, then for every $U\geq W$, we have
$$
x_U\in U\subseteq W.
$$
Now, using continuity of $\omega$ together with $\omega(e)=1$, both iterated limits $\Omega(x_U, x_V)$ converge to 1. By Proposition 2.1 in \cite{dl}, there exist subsequences $(x_{U_n})_{n\in{\Bbb N}}$ and $(x_{V_m})_{m\in{\Bbb N}}$ of $(x_U)_{U\in\frak{U}}\;$ such that
$$
\lim_n\lim_m\Omega(x_{U_n}, x_{V_m})=1=\lim_m\lim_n\Omega(x_{U_n}, x_{V_m}).
$$
Hence $\Omega$ can not be zero cluster, a contradiction. So zero clusters may exist only on discrete groups.$\hfill\square$\\

Let
$L^\infty (G,1/\omega)$ be the space of all measurable
functions $f$ on $G$ with
$$
\|f\| _{\infty,\;\omega}= \|f/\omega\|_\infty< \infty,
$$
where $\|.\|_\infty$ is the essential supremum norm. We denote by $L_0^\infty(G, 1/\omega)$ the
subspace of $L^\infty(G, 1/\omega)$ consisting of all
functions $f\in L^\infty (G,1/\omega)$ that vanish at
infinity.
It is well-known from \cite{mnr} that the dual space of
$L_0^\infty (G,1/\omega)^*$ is a Banach algebra with the first Arens product; see also \cite{lp, mmn}. One can show that $L_0^{\infty}(G, 1/\omega)^*$ is isomorphic with the set of all $F\in M_*(G, \omega)^*$ with
$$
\langle F, \lambda\rangle=\langle F, \lambda_0\rangle
$$
for all $\lambda\in M_*(G, \omega)^*$, where $\lambda_0=\lambda|_{L^1(G, \omega)}$; see \cite{mm}.

The first author and Moghimi \cite{mm} proved that if $M_*(G, \omega)^*$ is Arens regular, then $G$ is discrete. We are now in a position to prove the main result of this paper which is an improvement of Theorem 2 of \cite{rv} and is an answer to the open question presented in \cite{mm}.

\begin{theorem}\label{regular} Let $\omega$ be a weight function on a locally compact group $G$. Then the following assertions are equivalent.

\emph{(a)} $M_*(G, \omega)^*$ is Arens regular.

\emph{(b)} $L^1(G, \omega)$ is Arens regular.

\emph{(c)} $M(G, \omega)$ is Arens regular.

\emph{(d)} $L^1(G, \omega)^{**}$ is Arens regular.

\emph{(e)} $M(G, \omega)^{**}$ is Arens regular.

\emph{(f)} $L_0^\infty(G, \omega)^*$ is Arens regular.

\emph{(g)} $G$ is finite or $\Omega$ is zero cluster.

In this case, $G$ is discrete and countable.
\end{theorem}
{\it Proof.} Assume that $M_*(G, \omega)^*$ is Arens regular. Since $L^1(G, \omega)$ is a closed ideal in $M_*(G, \omega)^*$, it follows from Corollary 2.6.18 in \cite{d} that $L^1(G, \omega)$ is Arens regular. So (a) implies (b). It is well-known from \cite{rv} that $L^1(G, \omega)$ is Arens regular if and only if $G$ is finite or $G$ is discrete and $\Omega$ is zero cluster. From this and Lemma \ref{lem13} follows that the statements (b) and (g) are equivalent. From Theorem 4.7 and Corollary 4.11 in \cite{mnr} and Lemma \ref{lem13} we see that the other statements are equivalent.

To complete the proof, note that $G=\cup_{n=1}^\infty A_n$, where
$$
A_n=\{x\in G: \omega(x)\leq n\}.
$$
If $G$ is uncountable, then $A_m$ is infinite for some $m\in{\Bbb N}$. For every $x, y\in A_m$, we have $$\Omega (x, y)\geq 1/m^2.$$ This implies that $\Omega$ can not be zero cluster.$\hfill\square$

\begin{example}\label{ex}{\rm (i) Let $\alpha\geq 0$ and for every $n\in{\Bbb Z}$
$$
\omega_\alpha(n)= (1+|n|)^\alpha.
$$
One can prove that $\Omega_\alpha$ is zero cluster if and only if $\alpha>0$. It follows from Theorem \ref{regular} that $M_*({\Bbb Z}, \omega_\alpha)^*$ is Arens regular if and only if $\alpha>0$. In the case where $\alpha=0$, the Banach algebra $M_*({\Bbb Z})^*$ is not Arens regular.

(ii) Let $\alpha, \beta> 0$. For every $m, n\in{\Bbb Z}$ we define
$$
\omega(m, n)= (1+|m|)^\alpha(1+|n|)^\beta.
$$
Set $x_m=(m, 0)$ and $y_n=(0, n)$. Then $\Omega(x_m, x_n)=1$. Hence $M_*({\Bbb Z}^2, \omega)^*$ is not Arens regular.}
\end{example}

\begin{remark} {\rm Let $\omega$ be a weight function on a locally compact
group $G$. If $\Omega$ is either positive-cluster or $\Omega>\alpha$ for some $\alpha>0$, or $\omega$ is multiplicative, then $\Omega$ can not be zero cluster. So by Theorem \ref{regular}, the Banach algebra $M_*(G, \omega)^*$ is Arens regular if and only if $G$ is finite.}
\end{remark}

Baker and the second author \cite{br} gave the following result for the discrete convolution semigroup algebra $\ell^1 (S, \omega)$. In the following, we prove this result for locally compact groups.

\begin{corollary}\label{1} Let $\omega$ be a weight function on infinite locally compact
group $G$. Then the following assertions are equivalent.

\emph{(a)} $M_*(G, \omega)^*$ is Arens regular.

\emph{(b)} $L^1(G, \omega)$ is Arens regular.

\emph{(c)} $G$ is discrete and for every $A\subseteq G$ and each pair of sequences $(x_n)$ and $(y_n)$ in $G$,
$$
\{\chi_A (x_ny_m)\Omega(x_n, y_m): n<m\bar{\}}\cap\{\chi_A(x_ny_m)\Omega(x_n, y_m): n>m\bar{\}}\neq\emptyset.
$$

\emph{(d)} $G$ is discrete and for or each pair of sequences $(x_n)$ and $(y_n)$ in $G$ there exist subsequences $(a_n)$ and $(b_n)$ of $(x_n)$ and $(y_n)$, respectively, such that at least one of the following statements hold.

\quad\emph{(1)} $\lim_n\lim_m\Omega (a_n, b_m)=0=\lim_m\lim_n\Omega (a_n, b_m)$.

\quad\emph{(2)} either the rows or the columns of the matrix $(a_nb_m)$ are constant and distinct.

\quad\emph{(3)} The matrix $(a_nb_m)$ is constant.
\end{corollary}

As an immediate consequence of Corollary 3.5 in \cite{br} and Theorem \ref{regular}, we give the next result.

\begin{corollary}\label{alpha} Let $\omega_1$ and $\omega_2$ be weight functions on locally compact infinite group $G$. Then the following statements hold.

\emph{(i)} If $\Omega_1\geq\alpha\;\Omega_2$ for some $\alpha>0$ and $M_*(G, \omega_1)^*$ is Arens regular, then $M_*(G, \omega_2)^*$ is Arens regular.

\emph{(ii)} If there exist positive numbers $\alpha$ and $\beta$ such that $\alpha\omega_1\leq\omega_2\leq\beta\omega_1$, then $M_*(G, \omega_1)^*$ is Arens regular if and only if $M_*(G, \omega_2)^*$ is Arens regular.
\end{corollary}

Let us recall that a Banach algebra $\frak{A}$ is called a \emph{dual Banach algebra} if there exists a closed submodule $E$ of the dual module $\frak{A}^*$ such that $E^*=\frak{A}$. It is well-known that $(\frak{A}^{**}, \diamond)$ is a dual Banach algebra if and only if $\frak{A}$ is Arens regular; see for example Corollary 2.16 in \cite{dl}. From this together with Theorem \ref{regular}, we have the following result.

\begin{corollary} Let $\omega$ be a weight function on a locally compact
group $G$. Then $((M_*(G, \omega)^*)^{**}, \diamond)$ is a dual Banach algebra if and only if $G$ is finite or $\Omega$ is zero cluster.
\end{corollary}

For a weight function $\omega$ on $G$, we define $\omega^*(x)=\omega(x)\omega(x^{-1})$ for all $x\in G$. It is easy to see that $\omega^*$ is a weight function on $G$.

\begin{proposition}\label{co} Let $G$ be a locally compact group. Then the following assertions are equivalent.

\emph{(a)} $M_*(G)^*$ is Arens regular.

\emph{(b)} For every weight function $\omega$ on $G$, the Banach algebra $M_*(G, \omega)^*$ is Arens regular.

\emph{(c)} There exists a weight function $\omega$ on $G$ such that $M_*(G, \omega)^*$ is Arens regular and $\omega^*$ is bounded.

\emph{(d)}  There exists a weight function $\omega$ on $G$ such that $M_*(G, \omega)^*$ is reflexive and $\omega^*$ is bounded.

\emph{(e)} There exists a weight function $\omega$ on $G$ such that $M_*(G, \omega)^*$ is a $C^*-$algebra.

\emph{(f)} $G$ is finite.
\end{proposition}
{\it Proof.} First note that if $\omega$ is a weight function on $G$ such that $\omega^*$ is bounded, then there exists $\alpha>0$ such that
$$
\alpha\omega(x)\;\omega(y)\leq\omega(xy)
$$
for all $x, y\in G$. This shows that $\Omega$ can not be zero cluster. We also note that there exists a weight function $\omega$ on $G$ such that $\Omega$ can not be zero cluster. From These facts and Theorem \ref{regular}, we infer that  the assertions (a)-(d) and (f) are equivalent.  Now, let (e) hold. Then for every $x\in G$, we have
$$
\|\delta_x\ast\delta_x^*\|_\omega=\|\delta_x\|_\omega^2.
$$
This implies that $\omega=\Delta^{1/2}$, where $\Delta$ is the modular function of $G$. So
$$
\omega(xy)=\omega(x)\;\omega(y)
$$
for all $x, y\in G$. Therefore, $\Omega=1$. By Theorem \ref{regular}, (e) holds.$\hfill\square$\\

Let $\omega$ be a weight function on a locally compact
group $G$. By Proposition \ref{co}, if $M_*(G)^*$ is Arens regular,
then $M_*(G, \omega)^*$ is Arens regular. The converse,
however, is not true.

\begin{example}{\rm Let $\omega(n)= 1+|n|$ for all $n\in{\Bbb Z}$. Then $\Omega$ is zero cluster and so $M_*({\Bbb Z}, \omega)^*$ is Arens regular. But $M_*({\Bbb Z})^*$ isn't Arens regular}
\end{example}

As an immediate consequence of proposition \ref{co} we have the following result.

\begin{corollary}\label{sjmz} Let $\omega$ be a weight function on a locally compact
group $G$. Then the following assertions are equivalent.

\emph{(a)} $M_*(G, \omega)^*$ is Arens regular and $\omega^*$ is bounded.

\emph{(b)} $M_*(G, \omega)^*$ is reflexive and $\omega^*$ is bounded.

\emph{(c)} $G$ is finite.
\end{corollary}

Example \ref{ex}(i) shows that Corollary \ref{sjmz} is not true without the assumption that $\omega^*$ is bounded.

\section{\normalsize\bf Weight regularity of locally compact groups}

A locally compact group $G$ is called \emph{weight regular} if there exists a weight function $\omega: G\rightarrow [1,\infty)$ such that $M_*(G, \omega)^*$ is Arens regular.

\begin{theorem}\label{gr} Let $G$ be a locally compact
group. Then the following assertions are equivalent.

\emph{(a)} $G$ is weight regular.

\emph{(b)} $G$ is countable and discrete.


\emph{(c)} $G$ is finite or there exists a weight function $\omega: G\rightarrow [1, \infty)$ such that $\Omega$ is zero cluster.
\end{theorem}
{\it Proof.} The implications (a)$\Rightarrow$ (b) and (c)$\Rightarrow$(a) follow from Theorem \ref{regular}. The implication (b)$\Rightarrow$ (c) follows from Corollary 6.1.5 of \cite{dz}.$\hfill\square$\\

As a consequence of Theorem \ref{gr} we have the following result.

\begin{corollary}\label{z} Let $G$ be a locally compact infinite
group. If $G$ is compact or there exists a convergent net of distinct points of $G$, then $G$ is not weight regular. Furthermore, there is no weight function $\omega$ on $G$ such that $M_*(G, \omega)^*$ is Arens regular.
\end{corollary}
{\it Proof.} Let $(x_\alpha)_{\alpha\in A}$ be a convergent net of distinct points of $G$. If $G$ is weight regular, then $G$ is discrete. So $(x_\alpha)$ is eventually constant, a contradiction.
To complete, the proof note that if $G$ is an infinite compact group, then any net of distinct points of $G$, has a convergent subnet.
$\hfill\square$

\begin{example} {\rm By Theorem \ref{gr}, the additive group ${\Bbb Z}$ is weight regular, however, ${\Bbb R}$ and the tours group $${\Bbb T}=\{z\in{\Bbb C}: |z|=1\}$$ are not weight regular. So there is no weight function $\omega$ on ${\Bbb R}$ (respectively, ${\Bbb T}$) such that
$M_*({\Bbb R}, \omega)^*$ and $L^1({\Bbb R}, \omega)$ ( respectively, $M_*({\Bbb T}, \omega)^*$ and $L^1({\Bbb T}, \omega)$) are Arens regular.}
\end{example}

From Theorems \ref{regular} and \ref{gr} we have the following result due to Craw and Young \cite{cy}.

\begin{corollary} Let $G$ be a locally compact
group. Then there exists a weight function $\omega$ on $G$ such that $L^1(G, \omega)$ is Arens regular if and only if $G$ is countable and discrete.
\end{corollary}

\begin{proposition}\label{hom} Let $G_1$ and $G_2$ be locally compact groups  and $\psi: G_1\rightarrow G_2$ be a group homomorphism. Then the following statements are hold.

\emph{(i)} If $G_1$ is weight regular, then $\hbox{Im}\; \psi$ is weight regular.

\emph{(ii)} If $\psi$ is epimorphism and $G_1$ is weight regular, then $G_2$ is weight regular.

\emph{(iii)} If $\psi$ is epimorphism and $M_*(G_1)^*$ is Arens regular, then $G_2$ is weight regular.

\emph{(iv)} If $\psi$ is monomorphism and $M_*(G_2)^*$ is Arens regular, then $G_1$ is weight regular.

In these cases, $\psi$ is continuous.
\end{proposition}
{\it Proof.} Let $G_1$ be weight regular. Then there exists a weight function $\omega$ on $G_1$ such that $M_*(G_1, \omega)^*$ is Arens regular. Define the weight function $\omega_2$ on $\hbox{Im}\; \psi$ by
$$
\omega_2(\psi (t))= \hbox{inf}\;\omega_1(\psi^{-1}(\psi(t)))
$$
for all $t\in G_1$. Note that there is $0<\alpha<1$ such that for every $t\in G_1$
$$
\alpha\leq\omega_1(t)\leq\omega_2(\psi(t))+\alpha^2.
$$
So
$$
(1-\alpha)\omega_1(t)\leq \omega_2(\psi (t))\leq \omega_1(t).
$$
This implies that
$$
\Omega_2 (\psi (t), \psi (s))\leq\frac{1}{(1-\alpha ^2)}\Omega_1 (s, t)
$$
Now, Corollary \ref{alpha} proves (i). The statements (ii) and (iii) follow from (i).

Finally, let $M_*(G_2)^*$ be Arens regular. Then $G_2$ is finite. If $\psi$ is monomorphism, then $G_1$ is finite and so it is regular.$\hfill\square$\\

Let $\{G_i\}_{i\in I}$ be a family of locally compact groups and $\pi_j: \Pi_{i\in I}G_i\rightarrow G_j$ be the canonical projection, for $j\in I$. It is clear that $\pi_j$ is onto. Hence the following result holds.

\begin{corollary}\label{product} Let $\{G_i\}_{i\in I}$ be a family of locally compact groups. If $\Pi_{i\in I}G_i$ is  weight regular, then $G_i$ is weight regular for all $i\in I$.
\end{corollary}

Let us recall that a sequence $G_1\stackrel{f}\rightarrow G_2\stackrel{g}\rightarrow G_3$ of group homomorphisms is said to be \emph{exact} if $\hbox{Im}\; f=\hbox{ker}\; g$. An exact sequence of the form  $0\rightarrow G_1\stackrel{f}\rightarrow G_2\stackrel{g}\rightarrow G_3\rightarrow 0$ is called \emph{short exact}. If there exists a group homomorphism $h: G_2\rightarrow G_1$ such that $hf=1_{G_2}$, then the short exact sequence is called \emph{split}.

\begin{proposition}\label{soh} Let $G_1$, $G_2$ and $G_3$ be locally compact groups, $0\rightarrow G_1\stackrel{f}\rightarrow G_2\stackrel{g}\rightarrow G_3\rightarrow 0$ be a short exact sequence of group homomorphisms and $G_2$ be weight regular. Then the following statements hold.

\emph{(i)} $G_1$ is countable and $G_3$ is weight regular.

\emph{(ii)} If the given sequence is split, then $G_1$ and $G_3$ are weight regular.
\end{proposition}
{\it Proof.} Note that if the given sequence is short exact, then $g$ is onto. Also, if it is split, then $h$ is injective. These facts together with Theorem \ref{gr} and Proposition \ref{hom} prove the result.$\hfill\square$\\

In the sequel, we present a consequence of Proposition \ref{soh}.

\begin{corollary}\label{q} Let $N$ be a normal subgroup of locally compact
group $G$. Then the following statements hold.

\emph{(i)} If $G$ is  weight regular, then $G/N$ is weight regular and $N$ is is countable and open.

\emph{(ii)} If $G$ is weight regular and the sequence $0\rightarrow N\stackrel{\iota}\rightarrow G\stackrel{\pi}\rightarrow G/N\rightarrow 0$ is split, then $G/N$ and $N$ are weight regular, where $\iota$ is the inclusion map and $\pi$ is the quotient map.
\end{corollary}
{\it Proof.} It is easy to see that the sequence $0\rightarrow N\stackrel{\iota}\rightarrow G\stackrel{\pi}\rightarrow G/N\rightarrow 0$ is short exact. So $G/N$ is weight regular by Proposition \ref{soh}. From the weight regularity of $G$ and $G/N$ we infer that $G$ is countable and $G/N$ is discrete. Hence  $N$ is countable and open. So (i) holds. The statement (ii) follows at once from Proposition \ref{soh}.$\hfill\square$\\


We finish this section with the following result.

\begin{proposition} Let $G_i$ and $G^\prime_i$, for $i=1, 2, 3$, be locally compact groups and the sequences $0\rightarrow G_1\stackrel{f}\rightarrow G_2\stackrel{g}\rightarrow G_3\rightarrow 0$ and $0\rightarrow G_1^\prime\stackrel{f^\prime}\rightarrow G_2^\prime\stackrel{g^\prime}\rightarrow G_3^\prime\rightarrow 0$ be short exact. Let there exist group homomorphisms $\alpha_i: G_i\rightarrow G_i^\prime$, for $i=1, 2, 3$, such the obtained diagram
is commutative, i.e., $\alpha_2f=f^\prime\alpha_1$ and $\alpha_3g=g^\prime\alpha_2$. Then the following statements hold.

\emph{(i)} If $\alpha_1$ and $\alpha_3$ are group epimorphisms and $G_2$ is weight regular, then $G^\prime_2$, $G_3$ and $G^\prime_3$ are weight regular. Furthermore, $G_1$, $G_2$ and $G^\prime_1$ are countable.

\emph{(ii)} If $\alpha_1$ and $\alpha_3$ are group epimorphisms and $M_*(G_2)^*$ is Arens regular, then $G_i$ and $G^\prime_i$ are weight regular for $i=1, 2, 3$.

\emph{(iii)} If $\alpha_1$ and $\alpha_3$ are group monomorphisms and $M_*(G^\prime_2)^*$ is Arens regular, then $G_i$ and $G^\prime_i$ are weight regular for $i=1, 2, 3$.
\end{proposition}
{\it Proof.} (i) Since $g$ and $\alpha_3$ are surjective and the diagram is commutative, $g^\prime\alpha_2$ is surjective and hence
$$
\hbox{Im}\;g^\prime\alpha_2=C^\prime=\hbox{Im} g^\prime.
$$
So, if $b^\prime\in B^\prime$, then there exists $b\in B$ such that
$$
\alpha_2(b)- b^\prime\in\hbox{ker}\; g^\prime= \hbox{Im}\; f^\prime.
$$
But $\alpha_1$ is surjective and $f^\prime\alpha_1= \alpha_2 f$. Thus
$$
\hbox{Im}\; f^\prime=\hbox{Im}f^\prime\alpha_1=\hbox{Im}\;\alpha_2 f.
$$
Therefore,
$
\alpha_2(b)- b^\prime\in\hbox{Im}\;\alpha_2 f.
$
This shows that
$$
\alpha_2(b)- b^\prime=\alpha_2 f(a)
$$
for some $a\in A$. It follows that $\alpha_2(b- f(a))= b^\prime$. Hence $\alpha_2$ is surjective. Now, apply Proposition \ref{soh}.

(ii) This is an immediate consequence of (i).

(iii) By the assumption, $\alpha_3$ is injective and $
\alpha_3 g= g^\prime\alpha_2.
$
This implies that
$$
\hbox{ker}\;\alpha_2\subseteq \hbox{ker}\;g=\hbox{Im}\;f.
$$ By commutativity,
$
f^\prime\alpha_1= \alpha_2 f.
$
Since $f^\prime$ and $\alpha_1$ are injective, $\alpha_2 f$ is injective.
Hence $\alpha_2$ is injective. By Proposition \ref{soh}, the statement (iii) holds.$\hfill\square$



\section{\normalsize\bf Weighted function spaces}

Let $C_b(G)$ (respectively, $LUC(G)$) be the space of all bounded continuous (respectively, uniformly continuous) functions on $G$.
Let $C_b (G, 1/\omega)$ denote the space of all functions $f$ on $G$ such that $f/\omega\in C_b(G)$. A function $f\in C_b (G, 1/\omega)$ is called $\omega-$\emph{weakly almost periodic} (respectively, $\omega-$\emph{almost periodic}) if the set
$$
\{\frac{_xf}{\omega(x)\omega}: x\in G\}
$$
is relatively weakly (respectively, norm) compact in $C_b (G)$, where $_xf(y)=f(yx)$ for all $x, y\in G$. The set of all $\omega$(respectively, $\omega-$weakly) almost periodic on $G$ is denoted by $Ap(G, 1/\omega)$ (respectively, $Wap(G, 1/\omega)$). It is clear that
$$
Ap(G, 1/\omega)\subseteq Wap(G, 1/\omega)\subseteq C_b(G, 1/\omega).
$$
The equality may obtain for compact groups, however, it isn't necessary. Note that if $G$ is compact and $f\in C_b (G, 1/\omega)$, then the mapping $$x\mapsto\frac{_xf}{\omega (x) \omega}$$ from $G$ into $C_b (G)$ is continuous. This implies that $f\in Ap(G, 1/\omega)$. So the equality holds. In the sequel, we give necessary and sufficient condition under which the equality holds.

\begin{theorem}\label{wap} Let $\omega$ be a weight function on a locally compact
infinite group $G$. Then the following statements hold.

\emph{(i)} $Wap(G, 1/\omega)= C_b (G, 1/\omega)$ if and only if $G$ is compact or $\Omega$ is zero cluster.

\emph{(ii)} $Ap(G, 1/\omega)= C_b (G, 1/\omega)$ if and only if $G$ is either compact or discrete and $\Omega\in C_0(G\times G)$.
\end{theorem}
{\it Proof.} (i) Let $G$ be a non-compact group and $Wap(G, 1/\omega)= C_b (G, 1/\omega)$. Then
$$
LUC(G, 1/\omega)= C_b (G, 1/\omega),
$$
where $ LUC(G, 1/\omega)$ is the set of all $f\in C_b (G, 1/\omega)$ such that the map $x\mapsto\;_x(f/\omega)$ from $G$ into $C_b (G, 1/\omega)$ is norm continuous. Note that $Wap(G, 1/\omega)$ is a subspace of $LUC(G, 1/\omega)$. It is well-known from \cite{r} that $LUC(G, 1/\omega)= C_b (G, 1/\omega)$ if and only if $G$ is compact or discrete; see also \cite{b}. These facts show that $G$ is discrete. It follows from Corollary 3.8 (ii) in \cite{br} that $L^1(G, \omega)$ is Arens regular. By Theorem \ref{regular}, $\Omega$ is zero cluster.

Conversely, let $\Omega$ be zero cluster. In view of Lemma \ref{lem13} and Theorem \ref{regular}, $G$ is discrete and $L^1(G, \omega)$ is Arens regular. Applying Corollary 3.8 (ii) in \cite{br}, again, we have $Wap(G, 1/\omega)= C_b (G, 1/\omega)$.

(ii) Let $G$ be non-compact and $Ap(G, 1/\omega)= C_b (G, 1/\omega)$. Using (i) and Lemma \ref{lem13}, $G$ is discrete. Now, the result is proved if we only note that for discrete infinite group $G$,  $Ap(G, 1/\omega)= C_b (G, 1/\omega)$ if and only if $\Omega\in C_0(G\times G)$; see Corollary 3.18 (iii) in \cite{br}.$\hfill\square$


\begin{example}{\rm For every $n\in{\Bbb Z}$, we define $\omega(n)=1+|n|$. Then
$\Omega\in C_0({\Bbb Z}\times{\Bbb Z})$. So $$Ap({\Bbb Z}, 1/\omega)= C_b ({\Bbb Z}, 1/\omega)=Wap({\Bbb Z}, 1/\omega).$$
}
\end{example}

\begin{proposition}\label{mah3} Let $G$ be a non-compact group and $\omega$ be a weight function on $G$. Then the following assertions are equivalent.

\emph{(a)} $Wap(G, 1/\omega)= C_b (G, 1/\omega)$.

\emph{(b)} $M_*(G, \omega)^*$ is Arens regular.

\emph{(c)} $Wap(G, 1/\omega)= LUC(G, 1/\omega)$.

\emph{(d)} $\Omega$ is zero cluster.\\
\end{proposition}
{\it Proof.} It follows from Theorems \ref{regular} and \ref{wap} that the statements (a), (b) and (d) are equivalent. By Theorem 2.2 (ii) in \cite{r} and Lemma \ref{lem13} the statements (c) and (d) are equivalent. Finally, if $Ap(G, 1/\omega)= C_b (G, 1/\omega)$, then $Wap(G, 1/\omega)= C_b (G, 1/\omega)$. Hence $\Omega$ is zero cluster. By Theorem \ref{gr}, $G$ is weight regular.$\hfill\square$









\begin{theorem} Let $G$ and $G^\prime$ be non-compact groups, $N$ be a normal subgroup of $G$ and $\omega$ and $\omega_p$ be weight functions on $G$ and $G\times G^\prime$, respectively. Then the following statements hold.

\emph{(i)} If $\psi: G\rightarrow G^\prime$ is a group epimorphism and $Wap(G, 1/\omega)= C_b (G, 1/\omega)$, then there exists a weight function $\omega^\prime$ on $G^\prime$ such that $Wap(G^\prime, 1/\omega^\prime)= C_b (G^\prime, 1/\omega^\prime)$.

\emph{(ii)} If $Wap(G, 1/\omega)= C_b (G, 1/\omega)$, then there exists a weight function $\omega_q$ on $G/N$ such that $Wap(G/N, 1/\omega_q)= C_b (G/N, 1/\omega_q)$.

\emph{(iii)} If $Wap(G\times G^\prime, 1/\omega_p)= C_b (G\times G^\prime, 1/\omega_p)$, then there exist weight functions $\omega_0$ and $\omega_0^\prime$ on $G$ and $G^\prime$, respectively, such that $Wap(G, 1/\omega_0)= C_b (G, 1/\omega_0)$ and $Wap(G^\prime, 1/\omega_0^\prime)= C_b (G^\prime, 1/\omega_0^\prime)$.
\end{theorem}
{\it Proof.} Let $\psi: G\rightarrow G^\prime$ be a group epimorphism and $Wap(G, 1/\omega)= C_b (G, 1/\omega)$. It follows from Propositions \ref{mah3} that $G^\prime$ is weight regular. By Proposition \ref{hom}, $G^\prime$ is weight regular. Again, by Proposition \ref{mah3}, we obtain $Wap(G^\prime, 1/\omega^\prime)= C_b (G^\prime, 1/\omega^\prime)$ for some a weight function $\omega^\prime$ on $G^\prime$. So (i) holds. The statements (ii) and (iii) follow from Proposition \ref{mah3} together with Corollary \ref{q} and Corollary \ref{product}, respectively.$\hfill\square$\\

Let $\frak{A}$ be a Banach algebra. Then $f\in\frak{A}^*$ is called \emph{weakly almost periodic} (respectively, \emph{almost periodic}) if the map $a\mapsto af$ from $\frak{A}$ into $\frak{A}^*$ is weakly compact respectively compact, where $\langle af, b\rangle= \langle f, ba\rangle$ for all $b\in \frak{A}$. The spaces of all weakly almost periodic (respectively, almost periodic) functionals on $\frak{A}$ are denote by $WAP(\frak{A})$ and $AP(\frak{A})$, respectively.

\begin{theorem}\label{2122} Let $\omega$ be a weight function on a locally compact infinite group $G$. Then the following assertion are equivalent.

\emph{(a)} $WAP(M_*(G, \omega)^*)= (M_*(G, \omega))^{**}$.

\emph{(b)} $Wap(G, 1/\omega)= C_b (G, 1/\omega)$ and $G$ is discrete.

\emph{(c)} $\Omega$ is zero cluster.

\emph{(d)} $M_*(G, \omega)^*$ is Arens regular.
\end{theorem}
{\it Proof.} The implications (a)$\Leftrightarrow$ (d) and (c)$\Rightarrow$(d) follow from Theorem 2.14 in \cite{dl} and Theorem \ref{regular}. Let (d) hold. Since $G$ is infinite, by Theorem \ref{regular}, $\Omega$ is zero cluster. So (b) follows from Lemma \ref{lem13} and Theorem \ref{wap}(i). That is, (d)$\Rightarrow$(b). If (b) holds, then by Theorem \ref{wap}(i), $G$ is finite or $\Omega$ is zero cluster. By assumption, $\Omega$ is zero cluster. That is, (b)$\Rightarrow$(c).$\hfill\square$\\

\begin{theorem} Let $\omega$ be a weight function on a locally compact
infinite group $G$. Then the following assertion are equivalent.

\emph{(a)} $AP(M_*(G, \omega)^*)= (M_*(G, \omega))^{**}$.

\emph{(b)} $Ap(G, 1/\omega)= C_b (G, 1/\omega)$ and $G$ is discrete.

\emph{(c)} $G$ is discrete and $\Omega\in C_0(G\times G)$.
\end{theorem}
{\it Proof.} If $AP(M_*(G, \omega)^*)= (M_*(G, \omega))^{**}$, then $WAP(M_*(G, \omega)^*)= (M_*(G, \omega))^{**}$. By Theorem \ref{2122}, $G$ is discrete. Hence $AP(\ell^1(G, \omega))= \ell^\infty(G, 1/\omega)$. Therefore, $G$ is discrete and $\Omega\in C_0(G\times G)$. That is, (a) implies (c). By Theorem \ref{wap}, the statements (b) and (c) are equivalent. The implication (c)$\Rightarrow$(a) is clear.$\hfill\square$\\





Using a routine argument, the next result is established. So we omit it.

\begin{proposition} Let $\omega$ be a weight function on $G$. Then the following statements hold.

\emph{(i)} $C_0(G,1/\omega)= C_b (G, 1/\omega)$ if and only if $G$ is compact.

\emph{(ii)} $L^\infty(G,1/\omega)= C_b (G, 1/\omega)$ if and only if $G$ is discrete.

\emph{(iii)} $L^\infty_0(G,\omega)= C_b (G, 1/\omega)$ if and only if $G$ is finite.
\end{proposition}

\section{\normalsize\bf Amenability of $M_*(G, \omega)^*$}

Let us recall that the Banach algebra $M_*(G, \omega)^*$ is called \emph{amenable} if every continuous derivation from $M_*(G, \omega)^*$ into $E^*$ is inner for all Banach $M_*(G, \omega)^*-$module $E$.

\begin{theorem}\label{amenable} Let $\omega$ be a weight function on locally compact group $G$. Then the following assertions are equivalent.

\emph{(a)} $M_*(G, \omega)^*$ is amenable.

\emph{(b)}  $M(G, \omega)$ is amenable.


\emph{(c)} $G$ is a discrete amenable group and $\omega^*$ is bounded.

In this case, $M_*(G, \omega)^*= L^1(G)$.
\end{theorem}
{\it Proof.} Since $C_0(G, 1/\omega)$ is a closed subspace of $M_*(G, \omega)$, we imply that
$$
M_*(G, \omega)^*= M(G, \omega)\oplus M(G, \omega)_0,
$$
where
$$
M(G, \omega)_0=\{\Psi\in M_*(G, \omega)^*: \Psi|_{M(G, \omega)}=0\}.
$$
Let $\Phi\in M_*(G, \omega)^*$. Then there exists a net $(\mu_\alpha)_\alpha$ in $M(G, \omega)$ such that $\mu_\alpha\rightarrow\Phi$ in the weak$^*$-topology of $M_*(G, \omega)^*$. If $\Psi\in M(G, \omega)_0$ and $\nu\in M(G, \omega)$, then
$$
\langle\Phi\diamond\Psi, \nu\rangle=\lim_\alpha\langle\mu_\alpha\diamond\Psi, \nu\rangle=\lim_\alpha\langle\Psi, \nu\ast\mu_\alpha\rangle=0.
$$
Consequently, $M(G, \omega)_0$ is a left ideal in $M_*(G, \omega)^*$. On the hand, $\Phi=\mu+\Phi_0$ for some $\mu\in M(G, \omega)$ and $\Phi_0\in M(G, \omega)_0$. Then
$$
\Psi\diamond\Phi=\Psi\diamond(\mu+\Phi_0)=\Psi\diamond\mu+ \Psi\diamond\Phi_0.
$$
It is clear that $\Psi\diamond\mu\in M(G, \omega)_0$. Since $M(G, \omega)_0$ is a left ideal in $M_*(G, \omega)^*$, we have $\Psi\diamond\Phi_0\in M(G, \omega)_0$. So $M(G, \omega)_0$ is a right ideal in $M_*(G, \omega)^*$. Therefore, $M(G, \omega)_0$ is an ideal in $M_*(G, \omega)^*$.

Now, if $M_*(G, \omega)^*$ is amenable, then $$\frac{M_*(G, \omega)^*}{M(G, \omega)_0}\cong M(G, \omega)$$ is amenable; see for example \cite{r1}. So (a) implies (b).

The second author and Vishki \cite{rv} showed that $M(G, \omega)$ is amenable if and only if
$G$ is a discrete amenable and $\omega^*$ is bounded. Hence (b) and (c) are equivalent.

Let us recall that Gronback \cite{gro} proved that $L^1(G, \omega)$ is amenable if and only if $G$ is amenable and  $\omega^*$ is bounded. The first author and Moghimi \cite{mm} prove that $G$ is discrete if and only if $M_*(G, \omega)^*=L^1(G, \omega)$. These facts show that (c) implies (a).

It is well-known from \cite{w} that the mapping $\phi\mapsto \phi\;\omega$
 from $L^1(G, \omega)$ onto $L^1(G)$ is an isometric isomorphism of Banach spaces. Since $G$ is discrete, we have $$M_*(G, \omega)^*= L^1(G, \omega)=L^1(G).$$
So the assertions (a), (c) and (d) are equivalent.
$\hfill\square$

\begin{theorem} Let $\omega$ be a weight function on locally compact group $G$. Then the following assertions are equivalent.

\emph{(a)} $M_*(G, \omega)^*$ is amenable and Arens regular

\emph{(b)} $M_*(G, \omega)^{***}$ is amenable.

\emph{(c)} $M(G, \omega)^{**}$ is amenable.

\emph{(d)} $L^1(G, \omega)^{**}$ is amenable.

\emph{(e)} $G$ is finite.
\end{theorem}
{\it Proof.} From Theorems \ref{regular} and \ref{amenable} we infer that (a) and (e) are equivalent. Since $L^1(G, \omega)$ is an ideal in $M_*(G, \omega)^{*}$ and $M(G, \omega)$, respectively, it follows that $L^1(G, \omega)^{**}$ is an ideal in $M_*(G, \omega)^{***}$ and $M(G, \omega)^{**}$, respectively. Hence (b) and (c) imply (d). It follows from Theorem 4 in \cite{rv} that (d) implies (e). Trivially, (e) imply (b) and (c).$\hfill\square$

\vspace{2mm}

 {\footnotesize
\noindent {\bf Mohammad Javad Mehdipour}\\
Department of Mathematics,\\ Shiraz University of Technology,\\
Shiraz
71555-313, Iran\\ e-mail: mehdipour@sutech.ac.ir\\
{\bf Ali Rejali}\\
Department of Pure Mathematics,\\ Faulty of Mathematics and Statistics,\\ University of Isfahan,\\
Isfahan
81746-73441, Iran\\ e-mail: rejali@sci.ac.ir\\

\footnotesize

\begin{thebibliography}{99}

\bibitem{br} J. W. Baker and A. Rejali, On the Arens regularity of weighted convolution algebras, J. London Math. Soc., (2) 40 (1989) 535--546.

\bibitem{b} R. B. Burkel, Weakly Almost Periodic Functions on Semigroups, Gordon and Breach, New York, 1970.

\bibitem{cy} I. G. Craw and N. J. Young, Regularity of multiplication in weighted group and semigroup algebras, Quart. J. Math. Oxford, 25 (1974) 351--358.

\bibitem{dgh} H. G. Dales, F. Ghahramani and A. Y. A. Helemskii, The amenability of measure algebras, J. London Math. Soc., (2) 66 (2002) 213--226.

\bibitem{dl} H. G. Dales and A. T. Lau, The second duals of Beurling algebras, Mem. Amer. Math. Soc., 177 (836) (2005).

\bibitem{d} M. Daws, Arens regularity of the algebra of operators on a Banach space, Bull. London Math. Soc., 36 (2004) 493--503.

\bibitem{dsl} S. Degenfeld-Schonburg and R. Lasser, Multipliers on $L^p$-spaces for hypergroups, Rocky Mountain J. Math., 43 (4) (2013) 1115–-1139.

\bibitem{dh} J. Duncan and S. A. R. Hosseiniun, The second dual of a Banach algebra, Proc. Roy. Soc. Edinburgh, A 84 (1979) 309–-325.

\bibitem{dz} H. A. M. Dzinotyiweyi, The analogue of the group algebra for topological semigroups, Research Notes in Mathematics, 98. Pitman (Advanced Publishing Program), Boston, MA, 1984.

\bibitem{r1} H. R. Ebrahimi Vishki, B. Khodsiani and A. Rejali, Arens Regularity of certain weighted semigroup algebras and countability, Semigroup Froum, 92 (2016) 304--310.

 \bibitem{f} B. Forrest, Arens regularity and discrete groups, Pacific J. Math., 151 (2) (1991) 217--227.

 \bibitem{f1} B. Forrest, Arens regularity and the $A_p(G)$ algebras, Proc. Amer. Math. Soc., 119 (2) (1993) 595--598.

\bibitem{g1} C. C. Graham, Arens regularity and weak sequential completeness for quotients of the Fourier algebra, Illinois J. Math., 44 (4) (2000) 712--740.

\bibitem{g2} C. C. Graham, Arens regularity and the second dual of certain quotients of the Fourier algebra, Q. J. Math., 52 (1) (2001) 13--24.

\bibitem{g3} C. C. Graham, Arens regularity for quotients $A_p(E)$ of the Herz algebra, Bull. London Math. Soc., 34 (4) (2002) 457--468.

 \bibitem{gro} N. Gronback, Amenability of weighted convolution algebras on locally compact groups, Trans. Amer. Math. Soc., 319 (1990) 765--775.

\bibitem{lp} A. T. Lau and J. Pym, Concerning the second dual of the group
algebra of a  locally compact group, J. London Math. Soc., 41 (1990) 445--460.

\bibitem{lu} A. T. Lau and A. Ulger, Some geometric properties on the Fourier and Fourier-Stieltjes algebras of locally compact groups, Arens regularity and related problems, Trans. Amer. Math. Soc., 337 (1) (1993) 321--359.

\bibitem{mmn} S. Maghsoudi, M. J. Mehdipour and R. Nasr-Isfahani,
Compact right multipliers on a Banach algebra related to locally
compact semigroups, Semigroup Forum, 83 (2) (2011) 205--213.

\bibitem{mnr} S. Maghsoudi, R. Nasr-Isfahani and A. Rejali,
Strong Arens irregularity of Beurling algebras with a locally convex topology, Arch. Math., 86 (5) (2006) 437--448.

\bibitem{m} D. Malekzadeh Varnosfaderani, Derivations, Multiplers and Topological Centers of Certain
Banach Algebras Related to Locally Compact Groups, Ph.D. thesis,  University of Manitoba, 2017.

\bibitem{mm} M. J. Mehdipour and GH. R. Moghimi, The existence of non-zero compact right multipliers and Arens regularity of weighted Banach algebras, preprint.

\bibitem{r0} H. Reiter and J. D. Stegeman, Classical Harmonic Analysis and Locally Compact Groups, London Math. Society Monographs, 22, Clarendon Press, Oxford, 2000.

\bibitem{r} A. Rejali, Weighted function spaces on topological groups, Bull. Iranian Math. Soc., 22 (2) (1996) 43--63.

\bibitem{rv} A. Rejali and H. R. Vishki, Regularity and amenability of the second dual of weighted group algebras, Proyecciones, 26 (2007) 259--267.

\bibitem{r1} V. Runde, Lectures on amenability, Lecture Notes in Mathematics 1774, Springer Verlag, Berlin, 2002.

\bibitem{u1} A. Ulger, Arens regularity of the algebra $C(K,A)$, J. London Math. Soc., (2) 42 (1990) 354--364.

\bibitem{u2} A. Ulger, Some stability properties of Arens regular bilinear operators, Proc. Edinburgh Math. Soc., (2) 34 (1991) 443--454.

\bibitem{u} A. Ulger, Arens regularity of weakly sequentially complete Banach algebras, Proc. Amer. Math. Soc., 127 (1999) 3221--3227.

\bibitem{w} M. C. White, Characters on weighted amenable groups, Bull. London Math. Soc., 23 (1991) 375--380.

\end{thebibliography}
\end{document}